\documentclass[pdflatex,sn-mathphys-num]{sn-jnl}

\usepackage{graphicx}%
\usepackage{multirow}%
\usepackage{amsmath,amssymb,amsfonts}%
\usepackage{amsthm}%
\usepackage{mathrsfs}%
\usepackage[title]{appendix}%
\usepackage{xcolor}%
\usepackage{textcomp}%
\usepackage{manyfoot}%
\usepackage{booktabs}%
\usepackage[ruled]{algorithm2e}
\usepackage{listings}
\usepackage{comment}
\usepackage{enumitem}
\usepackage[short]{optidef}
\usepackage{comment}
\usepackage{subcaption}
\usepackage{tikz}
\usetikzlibrary{shapes.misc, positioning, shapes, fit}
\newtheorem{theorem}{Theorem}
\newtheorem{proposition}[theorem]{Proposition}%

\usepackage{background}

\begin{document}

\SetBgContents{Published at \url{https://doi.org/10.1287/ijoc.2024.0938}}      
\SetBgPosition{current page.center}
\SetBgAngle{0}                                    
\SetBgColor{gray}                                 
\SetBgScale{1.5}                                  
\SetBgHshift{0}                                   
\SetBgVshift{9cm} 

\title[A CG algorithm with DCA for MSSC]{A column generation algorithm with dynamic constraint aggregation for minimum sum-of-squares clustering}


\author[1]{\fnm{Antonio M.} \sur{Sudoso}}\email{antoniomaria.sudoso@uniroma1.it}

\author*[2]{\fnm{Daniel} \sur{Aloise}}\email{daniel.aloise@polymtl.ca}


\affil[1]{\orgdiv{Department of Computer, Control and Management Engineering ``Antonio Ruberti"}, \orgname{Sapienza University of Rome}, \orgaddress{\street{Via Ariosto 25}, \city{Rome}, \postcode{00185}, \country{Italy}}}

\affil[2]{\orgdiv{Department of Computer and Software Engineering}, \orgname{Polytechnique Montréal}, \orgaddress{\street{2500 Chem. de Polytechnique}, \city{Montréal}, \postcode{H3T 1J4}, \country{Canada}}}



\abstract{The minimum sum-of-squares clustering problem (MSSC), also known as $k$-means clustering, refers to the problem of partitioning $n$ data points into $k$ clusters, with the objective of minimizing the total sum of squared Euclidean distances between each point and the center of its assigned cluster.  We propose an efficient algorithm for solving large-scale MSSC instances, which combines column generation (CG) with dynamic constraint aggregation (DCA) to effectively reduce the number of constraints considered in the CG  master problem. DCA was originally conceived to reduce degeneracy 
in set partitioning problems
by utilizing an aggregated restricted master problem obtained from a partition of  the set partitioning  constraints into disjoint clusters. In this work, we explore the use of DCA within a CG algorithm for MSSC exact solution. Our method is fine-tuned by a series of ablation studies on DCA design choices, and is 
demonstrated to significantly outperform existing state-of-the-art exact approaches available in the literature.}

\keywords{Data clustering, Column generation, Dynamic constraint aggregation, Global optimization}



\maketitle

\section{Introduction \label{intro}}

The Minimum Sum-of-Squares Clustering (MSSC) problem is a fundamental problem in data analysis that has attracted significant attention from the fields of mathematical optimization and computer science throughout the years. The problem consists of partitioning a set of data points into \emph{clusters} such that the sum of squared distances from each point to the centroid of its cluster is minimized, which is equivalent to minimizing within-cluster variances.

A substantial amount of literature exists on this topic, with most of it centered around the classical and widely used ``$k$-means" algorithm \citep{forgy1965cluster}. The $k$-means heuristic, which has approximately 2 million references on Google Scholar in 2024, is an algorithm aimed at optimizing the MSSC problem. It starts with an initial partition of the data points into $k$ clusters. In the sequel, the centroids of the clusters are computed, and each point is assigned to its closest centroid. This process is then repeated until the partition is stable, i.e., no data point is reassigned to a different cluster centroid between two consecutive iterations of the algorithm. This iterative approach often leads to local optima solutions that can be far from the global optimum, particularly when dealing with a large number of clusters.  Additionally, the algorithm is highly sensitive to the initial partition.

The standard MSSC formulation is a {mixed-integer nonlinear programming} (MINLP) model, which represents the problem as the optimization of the discrete assignments of data points to continuous cluster centers to minimize the sum of squared distances between them.
Formally, the MSSC can be stated as follows:
\begin{subequations}
\begin{align}
\label{eq:MSSC}
\min \quad & \sum_{i=1}^n \sum_{c=1}^k x_{ic}\|p_i - y_c\|^2  \\
\textrm{s.\,t.} \quad & \sum_{c=1}^k x_{ic} = 1 \qquad \forall i \in \{1, \dots, n\}, \label{eq:MSSCa}\\
& x_{ic} \in \{0, 1\} \qquad \forall i \in \{1,\dots,n\}, \ \forall c \in \{1,\dots,k\},\\
& y_c \in \mathbb{R}^d \qquad \forall c \in \{1, \dots, k \}.
\end{align}
\end{subequations}
\noindent Here, $p_i \in \mathbb{R}^d$, for $i \in \{1,\dots,n\}$, where $d$ is the data dimension. The binary decision variable $x_{ic}$ expresses whether data point $i$ is assigned to cluster $c$ ($x_{ic}=1)$ or not ($x_{ic}=0$).
Constraints~\eqref{eq:MSSCa} ensure that each point is assigned to exactly one cluster. Although not explicitly expressed in the formulation, the centers of the $k$ clusters  $y_c$, for $c\in \{1,\dots,k\}$,  are located at the centroids of the clusters due to first-order optimality conditions on the $y$ variables.

Until very recently, a column generation (CG) algorithm published in 
\cite{aloise2012improved} was considered to be the state-of-the-art method for solving MSSC to global optimality. This shifted with the introduction of an SDP-based branch-and-bound algorithm proposed in~\cite{piccialli2022sos}, which was shown to solve real-world MSSC instances up to 4,000 data points. The main contributions of this work are:
\begin{itemize}
    \item We update the column generation of~\cite{aloise2012improved} in view of advances proposed in the literature for column generation algorithms in recent years,
more particularly with the use of Dynamic Constraint Aggregation (DCA)~\citep{elhallaoui2005dynamic,elhallaoui2010multi}.
    \item We conduct an extensive series of ablation experiments to explore various design choices in DCA for an optimization problem distinct from 
    the ones where DCA has been typically employed in the literature (i.e., vehicle routing problems \citep{elhallaoui2005dynamic}, crew scheduling \citep{desaulniers2020dynamic,elhallaoui2010multi,bouarab2017dynamic,saddoune2012integrated} and its variants). 
    \item Our results demonstrate that our CG algorithm enhanced by DCA outperforms current state-of-the-art exact methods on the tested problems, being  
    able to solve real-world instances with up to $\approx$  6000 data points.
   
    \item 
    Our CG+DCA implementation is entirely open-source and freely accessible,    
    with the exception of the commercial code utilized for solving linear relaxations. Previous implementations rely on the proprietary GENCOL column generation library.
   
\end{itemize}

We focused our experiments on clustering instances in the Euclidean plane, since the pricing problem of our CG algorithm can be very efficiently solved for $d=2$, as shown in~\cite{aloise2012improved}. This allows us to thoroughly evaluate the enhancements brought forth by DCA when applied to the solution of the CG restricted master problem. Clustering in two-dimensions has been used in various contexts and applications. For instance, it can be effectively employed in facility location problems to help identify optimal locations for facilities based on spatial data, representing customer locations or demand points. Once the clusters are formed, the centroids of these clusters can be considered as candidate locations for placing facilities. 
MSSC is in fact NP-hard for general $k$ in the plane in~\cite{mahajan2012planar}.

    Our paper is organized as follows. Section~\ref{sec:related} covers the related works regarding exact solution methods for MSSC.
    Section~\ref{sec:dca} {explains} the main ideas behind the DCA methodology, whereas Section~\ref{sec:dca-mssc}
    details our utilization of DCA within a CG algorithm  tailored for MSSC. Section~\ref{sec:bb} presents our branching strategy used to obtain integer clustering solutions and its impact on subproblem structure and solution. Section~\ref{sec:results} reports our computational experiments and results. Finally, final remarks with future research directions are provided in Section~\ref{sec:final}.

\section{Related works \label{sec:related}}

Given its relevance to several data mining applications, numerous optimization methods have been proposed in the literature to approach the MSSC problem, most of them by means of approximate methods. 
We contend, however, that global optimization is paramount in cluster analysis.  This is because clustering lacks the reliance on any form of side information to guide the formation of the clusters,
as opposed to supervised classification methods where a subset of data points are pre-classified to guide the underlying optimization process. As such, clustering analysis typically requires the interpretation of the obtained clusters by domain experts, which may be completely wrong shall they are presented to approximate clustering solutions of bad quality. This may lead to unexpected consequences in medical and financial contexts, for example.

The classes of methods for solving MSSC to global optimality differ on the formulations approached as well as on their bounding and branching strategies. 
To the best of our knowledge, the earliest use of branch-and-bound methods is attributed to \cite{koontz1975branch}, with subsequent refinement by \cite{diehr1985evaluation}. This procedure deals with partial solutions to the MSSC problem composed of fixed assignments for a subset $P' \subset P=\{p_1,\ldots,p_n\}$.  The fundamental observation used in this method is that the optimal solution value of the MSSC problem on $P$ is no less than the optimal solution of the MSSC problem on $P'$ plus the optimal solution value of the MSSC
on~$P - P'$. This approach was further pursued and improved by \cite{brusco2006repetitive}, who
developed the so-called repetitive-branch-and-bound algorithm (RBBA).
RBBA produces optimal solutions for synthetic datasets featuring up to 240 points and addresses challenging scenarios lacking inherent cluster structures, involving as many as 60 data points.

A more traditional line of research employs branch-and-bound algorithms where lower bounds are obtained through appropriate mathematical programming relaxations of the MINLP model in \eqref{eq:MSSC}. For instance, \cite{sherali2005global} proposed a branch-and-bound algorithm where lower bounds are determined by employing the reformulation-linearization-technique (RLT) \cite{sherali1998reformulation}. This technique involves transforming the underlying non-linear problem into an equivalent 0-1 mixed-integer program with a tight LP relaxation. 
More recently, \cite{burgard2023mixed} developed and tested several MIP techniques (e.g cutting planes, propagation techniques, branching rules, etc.) tailored to MSSC, so as to improve the performance of state-of-the-art MINLP solvers on solving \eqref{eq:MSSC} and its epigraph reformulation. Although numerical results show that the solution process of  SCIP solver is significantly improved, the authors were not able to match the performance of current state-of-the-art exact methods for the problem. 

A notable branch of literature focuses on employing semidefinite programming (SDP) techniques. \cite{peng2005new}, as well as  \cite{peng2007approximating}, established the equivalence between the MSSC problem and a 0-1 SDP reformulation. \cite{aloise2009branch} developed a branch-and-cut algorithm based on the linear programming (LP) relaxation of this 0-1 SDP model, achieving global optimality for instances with up to 202 data points.
More recently, \cite{piccialli2022sos} consider the SDP relaxation of the 0-1 SDP in \cite{peng2007approximating} and propose a new branch-and-bound algorithm that is capable of solving real-world instances with up to 4,000 data points. This algorithm features an efficient method for computing valid lower bounds from the SDP relaxation, which are strengthened by polyhedral cuts. The resulting relaxation is then solved in a cutting-plane scheme.
For the upper bound computation, the authors use the $k$-means algorithm with a new initialization procedure that exploits the solution of the SDP relaxation solved at each branch-and-bound node. 
To the best of our knowledge, this is the state-of-the-art solver for the MSSC problem.

Clustering  problems  can also be expressed as set partitioning problems by considering all possible clusters. Denote by ${T} = \{1, \dots, 2^{n}-1\}$ the index set of all possible subsets, i.e., clusters, of data points. For each $t \in {T}$, let $z_t$ be a binary variable indicating whether or not cluster $C_t$ is part of the solution, and let $a_{it}$ be a binary parameter that takes value 1 if the cluster $C_t$ contains the data point $p_i$, and 0 otherwise.
Denote by $y_t$ the centroid of points $p_i$ such that $a_{it} = 1$. 
The MSSC problem can then be rewritten by using the extended formulation:
\begin{subequations}
\label{prob:mp}
\begin{align}
\min \quad & \sum_{t \in {T}} c_t z_t  \\
\textrm{s.\,t.} \quad & \sum_{t \in {T}} a_{it} z_{t} = 1 \qquad \forall i \in \{1, \dots, n\}, \label{con:hard}\\
& \sum_{t \in {T}} z_{t} = k, \label{con:cardinality}\\
& z_t \in \{0, 1\} \qquad \forall t \in {T},
\end{align}
\end{subequations}
\noindent where $ c_t = \sum_{i=1}^n a_{it} \|p_i - y_t \|^2$.
Variable $z_t$ is equal to 1 if cluster $C_t$ belongs to the optimal partition, and to 0 otherwise. Constraints \eqref{con:hard} state that each entity belongs exactly to one cluster, and constraint \eqref{con:cardinality} expresses that the optimal partition contains exactly $k$ clusters. 
Without loss of generality, constraints \eqref{con:hard} 
can be replaced by
 $   \sum_{t \in {T}} a_{it} {z}_{t} \geq 1$, 
given that a covering of $P$, which is not a partition, cannot be an optimal MSSC solution.
Similarly, the constraint \eqref{con:cardinality} can be replaced by 
$\sum_{t \in {T}} {z}_{t} \leq k$,
since any MSSC partition with less than $k$ clusters has a cost greater or equal to the optimal MSSC partition with $k$ clusters.

 The extended formulation in \eqref{prob:mp} 
contains an exponential number of variables. Therefore, it cannot be solved in a straightforward way unless $n$ is small.
Formulation~\eqref{prob:mp} naturally lends itself to column generation (CG) procedures \citep{ lubbecke2005selected}. In this approach, the LP relaxation of the formulation is iteratively solved by employing a suitable \emph{pricing subproblem} to sequentially generate columns with negative reduced costs. These columns are then progressively appended to the \emph{restricted master program} (RMP), which  contains only a subset ${T'} \subseteq {T}$ of the columns from~\eqref{prob:mp}.
The~ \emph{primal}~RMP is  expressed as:
\begin{subequations}
\label{prob:rmp}
\begin{align}
\min \quad & \sum_{t \in {T'}} c_t z_t  \\
\textrm{s.\,t.} \quad & \sum_{t \in {T'}} a_{it} z_{t} \geq 1 \qquad \forall i \in \{1, \dots, n\}, \\
& \sum_{t \in {T'}} z_{t} \leq k, \\
& z_t \geq 0 \qquad \forall t \in {T'},
\end{align}
\end{subequations}
\noindent for which the reduced cost of a column $t \in T$ is given by 
 $    \pi_t = c_t - \sum_{i=1}^n a_{it}  \lambda_i + \sigma$, 
where the ${\lambda}_i$ for $i = 1, \dots, n$ and ${\sigma}$ are dual variables associated with the covering
constraints and with the side constraint, respectively.

After solving the linear relaxation of problem \eqref{prob:mp}, the solution obtained is examined to determine its integrality. As observed in~\cite{aloise2012improved}, the  relaxation of \eqref{prob:mp} is almost
always tight. However, if the solution is not an integer solution, branching becomes necessary.
The first branch-and-price algorithm for solving the MSSC extended formulation
was proposed by~\citet{du1999interior}, being later improved by~\citet{aloise2012improved} based on the observation that the pricing subproblem can be cast as a facility location problem with limited distances.
The authors proposed a series of valid inequalities based on geometric arguments that allowed to accelerate the solution of the pricing subproblem.

Nevertheless, large set partitioning problems, as the ones involved in data clustering, are difficult to solve by column generation due to the significant presence of degeneracy~\citep{elhallaoui2005dynamic}, that is, the average number of bases per extreme point of the feasible domain is very
large. This is even more pronounced for MSSC when the number of clusters $k$ is small.

\section{Dynamic Constraint Aggregation (DCA) \label{sec:dca}}

In an effort to reduce the impact of degeneracy while solving set partitioning problems by CG, \cite{elhallaoui2005dynamic}
proposed the so-called Dynamic Constraint Aggregation (DCA) technique. In DCA, the RMP contains significantly fewer constraints, expanding only when necessary to ensure optimality. This leads to a considerable decrease in the solution time for linear relaxations, along with less degenerate problems.

In DCA, the constraint set is represented by $W = \{w_1, \ldots,w_n\}$, where $w_i$ is the set partitioning constraint associated to the element $i$, and $\mathcal{Z}_i$ corresponds to the index set of columns that cover element $i$ in the current RMP. 
Utilizing these index sets $\mathcal{Z}_i$, we form a partition
$Q = \{G_1,\ldots,G_m\}$  of $W$ into $m < n$ groups  such that $G_j = \{ w \in W \mid \forall w_i, w_{i'} \in G_j : |G_j| > 1, \text{ it holds that } \mathcal{Z}_i = \mathcal{Z}_{i' } \}$, for $j=1,\ldots,m$. Finally, we denote by $\mathcal{I}_j$, for $j=1,\ldots,m$,  the  index set of the constraints in $G_j$. Note that, by consequence, $\mathcal{I}_j$ also refers to an index set of data points associated with those constraints. A CG algorithm using DCA  relies on an \emph{aggregated reduced master problem} (agRMP) where one \emph{aggregate constraint} $\bar{w}_j$ is used to replace the constraints in each set $G_j$, for $j=1,\ldots,m$.

The motivation of~\cite{elhallaoui2005dynamic} for proposing DCA stemmed from the observation that in applications such as vehicle routing problems, some sequences of deliveries are more likely to occur than others: the driver usually remains in the same vehicle during all trips of the day, i.e., many consecutive deliveries will remain grouped in the optimal solution. This grouping analysis allows for a considerable reduction in the size of the RMP in each iteration.
These consecutive tasks can then be seen as a single aggregated task, where the constraints corresponding to each task can be replaced by a single constraint.

As we mentioned in Section~\ref{sec:related}, the extended formulation of MSSC corresponds to a set partitioning problem where $n$ data points are partitioned into $k$ subsets called clusters. By leveraging the notation previously introduced in that section, we can express the agRMP of a CG for MSSC as: 
\begin{subequations}
\label{prob:agRMP}
\begin{align}
\min \quad & \sum_{t \in {T'}}  \bar{c}_t z_t  \\
\textrm{s.\,t.} \quad & \sum_{t \in T'} \bar{a}_{jt} z_{t} \geq 1 \qquad \forall j \in \{1, \dots, m\}, \label{con:hard_aggr}\\
& \sum_{t \in T'} z_{t} \leq k, \label{con:cardinality_aggr}\\
& z_t \geq 0 \qquad \forall t \in {T'},
\end{align}
\end{subequations}
\noindent where $\bar{a}_{jt}$ is a binary parameter that takes value 1 if the cluster $C_t$, associated to column $z_t$, contains all data points with indices in $\mathcal{I}_j$, and 0 otherwise. Similarly, the cost $\bar{c}_t$ of a column $z_t$ is computed taking into consideration all the data points in $\mathcal{I}_j$ for which $\bar{a}_{jt} = 1$.

The solution of the agRMP leads to the computation of an aggregated dual solution $\bar{\lambda} \in \mathbb{R}^m$, which subsequently requires disaggregation for generating negative reduced cost columns to be incorporated into the agRMP. These columns can be \emph{compatible} or \emph{incompatible} with the current partition $Q$. A  column $t \in T'$ is said compatible with $Q$ if, for all $j=1,\ldots,m$, it covers all elements in $\mathcal{I}_j$ or none of them, otherwise the column $t$ is said incompatible. 

When compatible negative reduced-cost columns are found, they can be pivoted into the basis of the current agRMP without modifying the aggregation, allowing the CG algorithm to advance to the subsequent iteration. However, if only incompatible negative reduced cost columns exist, an update of $Q$ becomes necessary. The CG algorithm concludes when no columns, whether compatible or incompatible with $Q$, are identified for addition as new columns to the current agRMP.  

Ideally, one would wish $m \ll n$ to ensure that the linear program remains small and degeneracy is kept low. However, we cannot know in advance which constraints could have been merged with respect to the optimal solution of~\eqref{prob:rmp}.

\section{CG algorithm with DCA  for MSSC \label{sec:dca-mssc}}

The intuition of DCA applied to the MSSC problem is that some data points have a higher likelihood of belonging to the same cluster. Specifically, data points that are very similar to each other (as determined by the Euclidean distance criterion) are more likely to be grouped together and can be treated as a single entity. Thus, the main idea is to group covering constraints related to similar data points and aggregating them into a single constraint.

In this section, we describe the main ingredients of our CG algorithm combined with DCA for the exact solution of MSSC, hereafter denoted CG+DCA.
Firstly, we detail how the initial partition of constraints $Q$ is obtained to start our algorithm. Secondly, we 
present the agRMP \emph{dual} problem and 
demonstrate how the dual variables associated to the covering constraints in the primal RMP can be bounded according to the partition of constraints 
$Q$, thus aiming to mitigate the unstable behavior of these variables throughout the algorithm. Thirdly, we describe how to disaggregate the aggregated dual variables after each CG iteration. 
 Fourthly, we express the pricing subproblem and the process of generating new columns for the agRMP, with a focus on maintaining their compatibility with
$Q$. Finally, we detail our approach for updating the partition 
$Q$ whenever incompatible columns need to be incorporated into the agRMP.

\subsection{Initial partition of covering constraints $Q$ \label{sec:initialQ}}

DCA relies on a partition of constraints $Q$ to be applied within CG.  
In the case of MSSC, the partition $Q$ is considered good if its subsets $G_j$, for $j=1,\ldots,m$ induce index sets $\mathcal{I}_j$ of data points that are part of the same column (i.e., cluster) in the optimal solution of \eqref{prob:rmp}.
 Thus, it appears natural to aggregate covering constraints associated to data points which are assigned to the same cluster on the basis of good 
  MSSC heuristic solutions. 
   In this work, a feasible (and supposedly good) MSSC solution $\bar{x} = \{C_1, \dots, C_k\}$ is obtained by executing the $k$-means heuristic. 

One straightforward initial partition $Q$ can be obtained from $\bar{x}$ by aggregating the covering constraints 
 into $m = k$ components, and thus, $G_j = \{ w_i \in W \mid p_i \in C_j\}$ for all $j \in \{1, \dots, k\}$. Nonetheless,  less aggressive initial constraint aggregations can also be explored. 
  To produce an initial  partition $Q$ with more than $k$ components, we proceed as follows. We first  execute  $k$-means for   
  $\tilde{k} > k$, thus obtaining the solution $\tilde{x} = \{\tilde{C}_1, \dots, \tilde{C}_{\tilde{k}}\}$. Then, upon examination of $\bar{x}$ and $\tilde{x}$, data points that are grouped together in the same cluster in both solutions have their corresponding covering constraints included as components of the initial $Q$. The remaining constraints represent singleton components of $Q$. 
  
 It is noteworthy that in this schema, the number of components $m$ of $Q$ can exceed $\tilde{k}$, as a cluster in  $\tilde{x}$ may contain data points belonging to multiple clusters in $\bar{x}$. 
In Section~\ref{sec:results}, we conduct an empirical investigation on the choice of $\tilde{k}$ for the CG+DCA algorithm.

\subsection{Bounds on aggregated dual variables \label{sec:bounds}}

A standard column generation algorithm is known to suffer from the  ``tailing off'' effect~\citep{lubbecke2005selected}. This refers to a phenomenon where the rate of improvement in the objective function value slows down as the algorithm progresses and more columns are added to the restricted master problem. This behavior is often observed in degenerate problems with a large number of variables, like the large-scale MSSC instances addressed in this paper. While there are various explanations for this phenomenon, a key factor behind it is the unstable behavior of the dual variables. 
The dual problem of the agRMP~\eqref{prob:agRMP} is given by:
\begin{subequations}
\label{prob:agRmp_dual}
\begin{align}
\max \quad & - k {\sigma} + \sum_{j=1}^{m} {\bar{\lambda}}_j  \\
\textrm{s.\,t.} \quad & -{\sigma} + \sum_{j=1}^m \bar{a}_{jt}{\bar{\lambda}}_j \leq \bar{c}_t \qquad \forall t \in {T'}, \\
& \sigma \geq 0, \quad \bar{\lambda}_j \geq 0 \qquad \forall j \in \{1, \dots, m\}, 
\end{align}
\end{subequations}
\noindent where $\bar{\lambda}_j$ represents the \emph{aggregated} dual variable associated to the aggregated constraint $\bar{w}_j$.

In order to address the tailing off issue mentioned at the beginning of this section, we consider a stabilization technique for the aggregated problem to regulate the behavior of the dual variables and ensure stability throughout the solution process. 
Let $\hat{l}_j \geq 0$ and  {$\hat{u}_j \geq \hat{l}_j$} be lower  and upper bounds, resp., on $\bar{\lambda}_j$
 for all $j \in \{1, \dots, m\}$.
 By using $\bar{\xi}_j$ to represent the dual variables associated with constraints
$\hat{l}_j \leq \bar{\lambda}_j$, and $\bar{\eta}_j$ to represent the dual variables associated with constraints 
$\bar{\lambda}_j \leq \hat{u}_j$, for $j=1,\ldots,m$, we obtain a stabilized version of the agRMP expressed as:
\begin{subequations}
\label{prob:sagRMP}
\begin{align}
\min \quad & \sum_{t \in T'} \bar{c}_t z_t - \sum_{j=1}^m \hat{l}_j \bar{\xi}_j + \sum_{j=1}^m \hat{u}_j \bar{\eta}_j  \\
\textrm{s.\,t.} \quad & -\bar{\xi}_j + \bar{\eta}_j + \sum_{t \in T'} \bar{a}_{jt} z_{t} \geq 1 \qquad \forall j \in \{1, \dots, m\}, \\
& \sum_{t \in T'} z_{t} \leq k, \\
& \bar{\xi}_{j}, \bar{\eta}_j \geq 0 \qquad \forall j \in \{1, \dots, m\},\\
& z_t \geq 0 \qquad \forall t \in T'.
\end{align}
\end{subequations}
\noindent The lower and upper bounds on the dual variables $\bar{\lambda}$ confine them within a box. If at the end of the column generation procedure, the dual optimum is attained on the boundary of the box, we have a direction towards which the box should be relocated. Otherwise, the optimum is attained in the interior of the box, producing the sought global optimum.

Estimated bounds for the dual variables associated
to non-aggregated constraints 
were presented in \cite{du1999interior}. Here, we adapt the estimation method to the presence of 
aggregated constraints.
The computation of these bounds is done in
the beginning of the column generation method with respect to the upper bound solution $\bar{x}$ 
and the initial partition $Q$ (see section~\ref{sec:initialQ}).
We describe below how we compute $\hat{l}_j$ and $\hat{u}_j$ for dual variables associated to aggregated constraints in $G_j \in Q$ such that $|G_j| > 1$.

To estimate the lower bound $\hat{l}_j$, we begin by assuming that only $\overline{\lambda}_j$ is upper bounded (i.e., $0 \leq  \overline{\lambda}_j \leq \hat{u}_j$ and $ \overline{\lambda}_{j'} \geq 0$ for $j' \in \{1,\ldots,m\} \setminus \{j\} $). Next, observe that if $\bar{\eta}_j >0$ this implies that $\overline{\lambda}_j \leq \hat{u}_j$ is active. Consequently, if $\hat{l}_j \leq \hat{u}_j$, then  $\hat{l}_j$ serves as a lower bound for $\overline{\lambda}_j$. Therefore, to determine a lower bound for $\overline{\lambda}_j$, we look for values of  $\hat{u}_j$ that ensure 
$\bar{\eta}_j >0$.

Let $v(\mathcal{I})$ be the value of the linear relaxation of~\eqref{prob:agRMP} when we are looking for a partition of $\mathcal{I} = \{\mathcal{I}_1,\ldots,\mathcal{I}_m\}$ into $k$ clusters. If $v(\mathcal{I}) - v(\mathcal{I}\setminus \{ \mathcal{I}_j\}) > \hat{u}_j$, then 
$\bar{\eta}_j$ will be strictly positive at the optimum of \eqref{prob:sagRMP}. This implies that
\[ \hat{l}_j < v(\mathcal{I}) - v(\mathcal{I}\setminus \{ \mathcal{I}_j\}). \]

To estimate $v(\mathcal{I}) - v(\mathcal{I}\setminus \{ \mathcal{I}_j\})$, we bound $v(\mathcal{I}\setminus \{ \mathcal{I}_j\})$ from above. Let $f(\overline{x})$ be the cost of the upper bound solution $\overline{x}$, and $c^1,\ldots,c^k$ the cost of its different clusters (i.e. $f(\overline{x}) = \sum_{i=1}^k c^i$).
Let $s$ be the index of the cluster in $\bar{x}$ that contain all data points with indices in $\mathcal{I}_j$,
and
$\underline{c}^{s}_j$ the MSSC cost of cluster $s$ when we omit the data points in $\mathcal{I}_j$. 
Thus, we have:
\[ v(\mathcal{I}\setminus \{ \mathcal{I}_j\}) \leq f(\overline{x}) - c^s + \underline{c}^{s}_j.
 \]

\noindent It follows that
\[ \hat{l}_j < v(\mathcal{I}) - 
f(\overline{x}) + c^s - \underline{c}^{s}_j \leq v(\mathcal{I}) - v(\mathcal{I}\setminus \{ \mathcal{I}_j\}).
\]

\noindent If we assume that no integrality gap exists and that $\overline{x}$ is optimal, then the difference $v(\mathcal{I}) - 
f(\overline{x})$ is equal to 0.
Therefore,
the lower bound  $\hat{l}_j$ of $\bar{\lambda}_j$ can be computed as $
\hat{l}_j = c^{s} - \underline{c}^{s}_j$.

Similarly for the upper bound $\hat{u}_j$, assuming that only $\overline{\lambda}_j$ is lower bounded, one can show that  
\[ \hat{u}_j > v(\mathcal{I} \cup \{\mathcal{I}_j\}) - v(\mathcal{I}). \]
\noindent Then noting that
\[ v(\mathcal{I}\cup \{ \mathcal{I}_j\}) \leq f(\overline{x}) + \overline{c}^{s'}_j- c^{s'},
 \]
where $\bar{c}^{s'}_j$ refers to the MSSC cost of a cluster $s' \neq s$ when all data points in $\mathcal{I}_j$ are added to it, we have that
\[ \hat{u}_j > f(\overline{x}) - v(\mathcal{I}) + \overline{c}^{s'}_j- c^{s'} \geq v(\mathcal{I} \cup \{\mathcal{I}_j\}) - v(\mathcal{I}). \]

\noindent Note that the tightest upper bound value for $\hat{u}_j$ is obtained by choosing $s'$ to minimize $\overline{c}^{s'}_j- c^{s'}$. Hence, by assuming no integrality gap and that $\overline{x}$ is optimal, the upper bound $\hat{u}_j$ of $\overline{\lambda}_j$ can be computed as $\hat{u}_j = \min_{s^ {\prime} \neq s} (\bar{c}^{s^{\prime}}_j - c^{s'})$.

The bounds on the dual variables $\bar{\lambda}_j$ just described are computed one
at a time for $j=1,\ldots,m$. They are still valid when all $\bar{l}_j$ and $\bar{u}_j$ are put together. 
They assume that no integrality gap exists and that $\bar{x}$ is optimal. Since that cannot be warranted in our solution process, 
 we must check for bounds which are active at the optimum,  modify them, and resume the CG algorithm if necessary. In our implementation we compute $\alpha_j = \hat{u}_j - \hat{l}_j$, for $j=1,\ldots,m$, and update $\hat{l}_j \leftarrow \max\{0,\hat{l}_j - \alpha_j/2\}$ and 
 $\hat{u}_j \leftarrow \max\{\hat{u}_j + \alpha_j/2\}$.

\subsection{Dual variable disaggregation \label{sec:disaggregation}}

One challenge associated with the DCA method is that it does not provide a complete dual solution $(\lambda_1, \dots, \lambda_n)$ when solving the aggregated RMP. Instead, it offers an aggregated dual solution $\bar{\lambda}_j$ for each $G_j \in Q$, with $j = 1, \dots m$, that must be disaggregated to be used as input to the pricing subproblem.
The choice of the dual variable disaggregation procedure allows for considerable flexibility, with the only requirement being the complementarity between primal and dual solutions. For the MSSC problem, a sufficient condition to ensure that complementarity is met is the following:
\begin{equation}\label{constr:sum_compl}
    \sum_{i \in \mathcal{I}_j} \lambda_i = \bar{\lambda}_j \quad \forall j \in \{1, \dots, m\}.
\end{equation}


Indeed, one can employ various approaches to find a feasible solution for the linear system presented in \eqref{constr:sum_compl}, which possesses an infinite number of feasible solutions. \cite{elhallaoui2005dynamic} suggested that instead of finding a feasible solution of the system in \eqref{constr:sum_compl} it is preferable to add some inequalities enforcing that the reduced costs of a large known subset of incompatible columns with respect to $Q$ are non-negative.
For instance, \cite{elhallaoui2005dynamic, elhallaoui2010multi} applied DCA to a CG for the vehicle and crew scheduling problems where they assume that the tasks (set partitioning constraints) are ordered for their starting time. After performing a variable substitution, the authors showed that the linear system corresponds to the optimality conditions of a shortest path problem. In \cite{bouarab2017dynamic}, a linear problem rather than a shortest path problem is used to disaggregate the dual variables. This LP is known as the complementary problem and it is solved at every CG iteration to generate more central disaggregate dual variables.

An attractive disaggregation strategy is the sparse one, in which the dual value related to a single constraint in $G_j$ is set to $\bar{\lambda}_j$, while the remaining dual values related to the other constraints in $G_j$ are set to zero. 
Indeed, this disaggregation strategy reduces considerably the number of variables in the pricing problem (see Section~\ref{sec:pricing}).
Yet, another simple disaggregation strategy consists in equally distributing the value of $\bar{\lambda}_j$ among the dual values associated to the constraints in $G_j$, i.e.
\begin{equation}
    \label{eq:disagg_avg}
        \lambda_i = \frac{1}{|\mathcal{I}_j|} \bar{\lambda}_j, \quad \forall j \in \{1, \dots, m\}, \ i \in \mathcal{I}_j.
\end{equation}

\subsection{Pricing problem \label{sec:pricing}}

Once a complete dual solution $\lambda$ is obtained by means of a disaggregation procedure, CG+DCA can proceed by looking for new negative reduced columns to be added to the aggregated RMP. For that purpose, a pricing problem must be solved.
As mentioned in Section~\ref{sec:related}, the reduced cost of a column $t \in {T}$ is given by
 $   \pi_t = c_t + \sigma - \sum_{i=1}^n a_{it} \lambda_i$.
Note that this equation remains the same whether the RMP is stabilized or not (see Section~\ref{sec:bounds}).

Since we are interested in finding violated constraints $\pi_t < 0$, the pricing problem can be formulated as
\begin{equation}
    \pi^\star = \sigma + \min_{y_v \in \mathbb{R}^d, v \in \{0, 1\}^n} \sum_{i=1}^n \left(\| p_i - y_v \|^2 - \lambda_i  \right) v_i, \label{prob:pricing}
\end{equation}
where $y_v$ denotes the centroid of points $p_i$ for which $v_i = 1$. Note that problem~\eqref{prob:pricing} is cubic, and hence a nonconvex MINLP, since the objective function uses multiplications of variables $v_i$ with the norms that depend on the centroids $y$, which are variables of the problem as well. However, it can also be re-written as a convex MINLP in a lifted space by using its epigraph formulation \citep{aloise2012improved}.

Several methods have been employed in the literature to solve problem~\eqref{prob:pricing}. By leveraging the fact that the sum of squared distances from each data point within a given cluster to its centroid equals the sum of squared distances between every pair of data points within that cluster divided by its cardinality,
\citet{du1999interior} expressed  problem~\eqref{prob:pricing} as a hyperbolic (fractional) program in 0-1 variables with quadratic numerator and linear denominator as:

\begin{equation}\label{eq:aux:hyperbolic}
 \sigma + \min_{v \in \{0, 1\}^n}  \frac{\displaystyle \sum_{i=1}^{n-1} \sum_{j=i+1}^n \left(\|p_i - p_j\|^2 -\lambda_i - \lambda_j \right) v_{i} v_{j} - \sum_{i=1}^n \lambda_i v_i}{\displaystyle \sum_{i=1}^n v_i}, 
\end{equation}

\noindent which is solvable by a Dinkelbach-like algorithm that relies on a branch-and-bound algorithm
for unconstrained quadratic 0-1 optimization problems (see~\cite{du1999interior} ).
To accelerate the resolution of these problems, a Variable Neighborhood Search (VNS) heuristic~\citep{hansen2017variable} was used instead  as long as it was able to a find negative reduced cost columns.

In~\cite{aloise2012improved}, the pricing problem is solved by an approach based on geometrical arguments. We note that 
problem \eqref{prob:pricing} minimizes the sum of squared distances between the cluster center $y_v$ and each entity, subject to constraints on each distance. If $\|p_i - y_v\|^2 \leq \lambda_i$, then $v_i$ is equal to 1 in the optimum, and 0 otherwise. In the plane, this is equivalent to the condition that $v_i = 1$ if $y_v$ belongs to a disc $D_i = \{y \in \mathbb{R}^2 : \|p_i - y \| \leq \sqrt{\lambda_i} \}$, (i.e., a disc with radius $\sqrt{\lambda_i}$ centered at $p_i$), and 0 otherwise. 
By focusing on components of $v$ that are equal to 1, we then consider subproblems of the following type:
\begin{equation}
\label{prob:y}
    \tau(S) = \min_{y \in \mathbb{R}^2}  \sum_{i \in S} \| p_i - y \|^2  \quad \textrm{s.t.} \quad \|p_i - y \|^2 \leq \lambda_i \quad \forall i \in S ,
\end{equation}
where $S \subseteq \{1, \dots, n\}$ is a non-empty set.
Subproblems of type \eqref{prob:y} are convex programming problems that can be solved in closed form and whose optimal solution can be obtained by taking the barycenter of all the data points satisfying the constraint $\|p_i - y \|^2 \leq \lambda_i$ for all $i \in S$.

In \cite{aloise2012improved}, the authors demonstrate the existence of a set $S^{\ast} \subseteq \{1,\ldots,n\}$ for which the optimal $y^{\ast}$ of $\tau(S^{\ast})$ is likewise optimal for \eqref{prob:pricing}. In the two-dimensional Euclidean space, this implies that searching for the optimal solution  $(v^{\ast},y_v^{\ast})$ for~\eqref{prob:pricing}  can be reduced to solving subproblems of type~\eqref{prob:y} whose sets $S$ are formed by the indices of discs $D_i$ that intersect in the plane, as well as singleton sets $S$.  
As demonstrated by~\cite{venkateshan2020note}, the number of such subproblems is polynomially bounded by $O(n^2)$.

It is important to note that, when employing DCA within a CG algorithm, the pricing problem must primarily return compatible columns with negative reduced cost. While this restriction could alter the definition of our pricing problem,  it can be handled without modifying the structure of the problem. More precisely, after solving a subproblem $\tau(S)$, it is sufficient to verify whether $S$ is compatible with the aggregating partition $Q$. Nevertheless, it might happen that such a ``compatible'' set $S$ is no longer found, even though incompatible columns with negative reduced cost exist. In that case, $Q$ and the aggregated RMP must be updated to accommodate new incompatible columns.

\subsection{Aggregating partition update \label{sec:q_update}}

To ensure the exactness of the method, the aggregating partition $Q$ must be dynamically updated, and thus, a new aggregated restricted master problem must be solved. To this end, we now describe how the current aggregating partition can be disaggregated. 

At each CG iteration, the pricing problem returns a set of negative reduced cost columns $T_P$ and a set of columns $C_Q \subseteq T_P$ that are compatible with $Q$. When $C_Q = \emptyset$ and $T_P \neq \emptyset$, it means that none of the columns in $T_P$ are compatible. Therefore, we must select (at least) one of these columns to update $Q$, enabling the pricing algorithm to generate compatible columns once again. 

According to \cite{elhallaoui2010multi}, incompatible columns can be ranked according to their number of incompatibilities, i.e., the number of additional components that would need to be created in order to make these columns compatible with the current constraint aggregation. Formally, given an aggregating partition ${Q}$ and a column $t \in T$, the number of incompatibilities of $t$ with respect to ${Q}$ is given by
$    u(t, {Q}) = \sum_{j=1}^{m} \chi_{\Omega}(t, \mathcal{I}_j)$,
where $\chi_{\Omega}$ is the characteristic function of the set $\Omega = \{(t, \mathcal{I}) : t \cap \mathcal{I} \neq \emptyset \ \textrm{and} \ t \cap \mathcal{I} \neq \mathcal{I} \}$. That is, $\chi_{\Omega}(t, \mathcal{I}_j)$ is equal to 1 if the column $t$ contains some of the data points with indices in $\mathcal{I}_j$, but not all of them, and 0 otherwise. A column $t \in T$ and its associated variable $z_t$ are said to be $p$-incompatible with respect to ${Q}$ if $u(t, {Q}) = p$. One can easily verify that compatible columns are qualified as 0-incompatible columns.

To pivot an incompatible column into the basis of the agRMP, we update the aggregation by increasing the number of components (via additional covering constraints) and, thus, the size of this problem. 
However, as reported in \cite{elhallaoui2010multi}, if the choice of the incompatible column used for disaggregating $Q$ is made without considering its impact on the size of the resulting aggregated problem, the problem size may increase rapidly. 

In the computational experiments outlined in Section~\ref{sec:results}, we assess two methodologies for updating $Q$ when the pricing problem can no longer generate compatible columns, albeit still returning incompatible ones. One considers the minimum reduced cost incompatible column and another the incompatible column with the smallest number of incompatibilities.
 When multiple columns have the same number of incompatibilities, we choose the one with the smallest reduced cost. 
 The pseudo-code of the algorithm for updating the aggregating partition $Q$ from an incompatible selected column is presented in Online Appendix A.

\subsection{Algorithm overview}

The CG+DCA for MSSC is shown in Algorithm \ref{alg:dcapseudocode}. It begins by choosing an initial partition ${Q}$ (Section~\ref{sec:initialQ}).
 This solution is also used to estimate lower and upper bounds for the dual variables and construct the stabilized agRMP (Section~\ref{sec:bounds}). Then, it performs two types of iterations, with minor iterations being contained within major ones. A minor iteration consists of three steps. The current agRMP is solved using an LP solver to compute a primal solution and a dual (aggregated) solution. Here, the aggregated dual solution contains a dual variable for each aggregated  covering constraint in the agRMP. Next, the dual variable values are disaggregated (Section~\ref{sec:disaggregation}).
 Given these disaggregated dual values, the pricing problem is solved (Section~\ref{sec:pricing}) to generate negative reduced cost columns that are compatible with $Q$. There are three possible outcomes at this point. The first is when no negative reduced cost columns are produced, implying that the algorithm has reached an optimal solution for the original relaxed problem. The second is when compatible columns are found, in which case they are added to the agRMP, and a new minor iteration starts. Finally, if no compatible columns with negative reduced cost are generated, partition ${Q}$ must be modified (Section~\ref{sec:q_update}) to allow for some incompatible columns to become compatible before starting a new major iteration. When the CG algorithm converges, the optimal dual values are checked to see whether they are attained at the boundary of the estimated bounds. If so, some bounds are updated and the CG+DCA algorithm is resumed.

\begin{algorithm}[!hbt]
\footnotesize
\caption{CG+DCA for MSSC}
\label{alg:dcapseudocode}
\textbf{Input}: Data points $P$, number of clusters $k$.
\begin{enumerate}[label*=\arabic*., nolistsep, rightmargin=1cm] 
    \item Initialize the aggregating partition $Q$.
    \item Estimate lower and upper bounds for $\bar{\lambda}_1, \dots, \bar{\lambda}_m$, and build the agRMP with respect to $Q$.
    \item Solve agRMP and let $\textrm{MSSC}^\star$ be the optimal value.
    \item Disaggregate dual variables $\bar{\lambda}_1, \dots, \bar{\lambda}_m$ and obtain $\lambda_1, \dots, \lambda_n$.
    \item Generate a set of negative reduced cost columns $T_P$ and a subset $C_Q \subseteq T_P$ of compatible columns by solving the pricing problem.
    \begin{enumerate}
        \item If $T_P = \emptyset$ then Stop.
        \item If $C_Q = \emptyset$ then update $Q$ and go to Step 2.
        \item If $C_Q \neq \emptyset$ then add the columns in $C_Q$ to RMP and go to Step 3.
    \end{enumerate} 
    \item  If dual optimum is attained on the boundary of the box then update the box and go to Step 3. Otherwise Stop.
\end{enumerate}
\textbf{Output}: Valid lower bound $\textrm{MSSC}^\star$.
\end{algorithm}

\section{Branching}
\label{sec:bb}

To obtain integer solutions, the CG+DCA algorithm is embedded within a branch-and-bound framework. This process becomes an integral part of the branch-and-price framework, where branching divides the problem into subproblems that are handled by CG. We apply the Ryan-Foster branching whenever the optimal solution of the master problem is not integer \citep{du1999interior, aloise2012improved}. It consists of finding two rows $i_1$, $i_2$ such that there are two columns $t_1$ and $t_2$ with fractional values at the optimum and such that $a_{i_1 t_1} = a_{i_2 t_1} = 1$ and $a_{i_1 t_2} = 1$, $a_{i_2 t_2} = 0$. Such rows and columns necessarily exist in any fractional optimal solution of the continuous relaxation. 
Then, pairwise constraints are introduced in the pricing problem of both branching subproblems in the form of a must-link constraint for one branch, and as a cannot-link constraint in the other.
Hence, given the dual variables $(\sigma, \lambda_1, \dots, \lambda_n)$, the pricing problem in the presence of branching constraints can be expressed as
\begin{align}
\label{prob:branching_pricing}
\begin{split}
\min_{y_v \in \mathbb{R}^d, v \in \{0, 1\}^n} \quad & \sigma + \sum_{i=1}^n \left(\|p_i - y \|^2 - \lambda_i \right) v_i \\
\textrm{s.\,t.} \quad & v_i - v_j = 0 \qquad \forall (i, j) \in \textrm{ML}, \\
& v_i + v_j \leq 1 \qquad \forall (i, j) \in \textrm{CL},
\end{split}
\end{align}
\noindent where $\textrm{ML}$ and $\textrm{CL}$ are the index sets of pairs of data points involved in must-link and cannot-link constraints, respectively.
Similarly to~\eqref{prob:pricing}, problem~\eqref{prob:branching_pricing} is a nonconvex MINLP that can be approached in several different ways as detailed in Section~\ref{sec:pricing}. 
A VNS heuristic is used in  \cite{du1999interior,aloise2012improved} to speed up the process of finding negative reduced cost columns until optimality of~\eqref{prob:branching_pricing} must be checked by an exact optimization algorithm. 
 
We  introduce here a new heuristic  algorithm for~\eqref{prob:branching_pricing} that exploits its structure. In fact, the structure of the objective function and constraints imply that the computation of the global minimum with respect to each block of variables (i.e., binary variables $v_i$ and continuous variables $y$) can be done efficiently. More in detail, if we fix $v_i$ then the problem becomes an unconstrained optimization problem whose optimal solution can be obtained in closed form. On the other hand, if we fix $y$ to $\bar{y}$ and optimize with respect to $v_i$,  problem~\eqref{prob:branching_pricing} becomes an integer programming (IP) problem. 
Under mild assumptions, the next proposition states that 
this problem
can be solved in polynomial time. Its proof is provided in Online Appendix B.
\vspace{0.2cm}
\begin{proposition}
\label{prop:tum}
    Let $G = (V, E)$ be an undirected graph   
    with one node for each data point appearing in the constraint set CL, and with an edge between two nodes if the corresponding points appear together in a cannot-link constraint. Denote by $A \in \mathbb{R}^{(|CL|+|ML|) \times n}$ the constraint matrix of Problem \eqref{prob:branching_pricing}. If $G$ admits a bipartition then $A$ is totally unimodular.
\end{proposition}
\vspace{0.2cm}

Our heuristic pricing algorithm for the branching nodes works by looping over the following steps until convergence is reached:
\begin{enumerate}
    \item  Solve 
         \eqref{prob:branching_pricing} for fixed $\bar{y}$. Let $\bar{v}_i$ be the optimal solution and $Z = \{p_i : \bar{v}_i = 1\}$.
        \item Update $\bar{y}$ as the centroid of the data points in $Z$.
    \end{enumerate}

The initial guess for  $\bar{y}$ is obtained by using one column $t \in T_P$ returned by the unconstrained pricing (see Section~\ref{sec:pricing}). Consequently, the procedure is repeated for each negative reduced cost column in $T_P$. 
If the obtained column, formed by the members of $Z$, has negative reduced cost and is compatible with $Q$, it is added to the agRMP.

The convergence of our heuristic can be proved by adapting the convergence proof of the $k$-means algorithm. Specifically, convergence is guaranteed because the variables $v_i$ can only assume a finite number of values due to the finite number of possible assignments. Consequently, as the algorithm finds the global minimum with respect to the centroid $y$ for a given assignment $v_i$, and vice versa, there will be a point in the iteration sequence where the centroid and assignment variables no longer change. 

Our heuristic method is very efficient particularly when problem \eqref{prob:branching_pricing} for fixed $\bar{y}$ is solvable by linear relaxation. However,  this may not be always true along the enumeration tree.
 In fact, as more cannot-link constraints are added within a given branch, it becomes increasingly difficult for $G$ to admit a bipartition, as stated in Proposition \ref{prop:tum}. Fortunately, as we will observe in Section~\ref{sec:results}, most of the MSSC problems can be directly solved at the root node or with a few nodes by our branch-and-price algorithm.

Finally, if our heuristic fails to identify negative reduced cost columns compatible with the current aggregation 
$Q$, it might still identify  incompatible negative reduced cost columns to be added to the restricted master problem at the branching node. In such a scenario, 
$Q$ is updated as described in Section \ref{sec:q_update}, and the CG algorithm resumes. If even incompatible negative reduced cost columns are unavailable, problem \eqref{prob:branching_pricing} must be solved to global optimality. In this case, we follow the decomposition approach outlined in \cite[pp.~205-207]{aloise2012improved}, where subproblems are solved by means of a Dinkelbach-like algorithm.

\section{Computational results}\label{sec:results}

Our CG+DCA algorithm is implemented in C++. At each CG iteration, the agRMP is solved with the simplex algorithm implemented in Gurobi (version 11.0). We also use Gurobi to solve the LPs within the pricing heuristic and the binary QPs within the Dinkelback algorithm when exact pricing is needed in the branching nodes.
Initial upper bound solutions (UB) are obtained by the $k$-means algorithms executed 10,000 times from different seeds.
 The best-obtained clustering solution is used to build the initial set of $k$ columns, estimate box constraints for the dual variables, and initialize the aggregating partition. We run the experiments on a machine equipped with Intel(R) i7-13700H CPU and 32 GB of RAM under Ubuntu 22.04.3. To promote reproducibility, our implementation is open-source and freely accessible at \url{https://github.com/INFORMSJoC/2024.0938} \citep{SudosoAloise2025}.

We use several data instances from the TSPLIB
~\citep{reinelt1991tsplib}.
They consist of sets of towns with 2D coordinates, where each town correspond to a data point to be clustered.  The number of data points can be easily guessed from the name of the data instances. For example, \texttt{ch150} corresponds to an instance with $n=150$ data points. It is important to notice that, in principle, the used TSP data sets do not exhibit hidden cluster structures. As pointed out in~\cite{aloise2011evaluating}, the hardness of an MSSC instance for exact methods is not directly related to the number of data points, but rather to the separation of the clusters to be found.

The experiments reported in this section serve three main purposes. First, we aim to show that the incorporation of DCA in the CG method of~\cite{aloise2012improved} allows to reduce the the time spent in solving its RMP, thus surpassing the performance of the referenced method. Secondly, we explore through a series of ablation studies the possible design choices in the DCA methodology, analyzing their impact in the performance of the CG+DCA method introduced here. Finally, we conclude the experiments by evaluating the new exact method in comparison with current state-of-the-art solver for MSSC.

\subsection{DCA evaluation}

In this section, we evaluate the effectiveness of the DCA within a CG algorithm for MSSC. To this end, we can restrict our evaluation to the solution of the root node only without prejudice.

Tables~1-4 in Online Appendix C report the results obtained by our CG+DCA algorithm at the root node for MSSC under different aggregation levels (column $\ell$), which are identified as follows:
\begin{itemize}
\item $\ell = n$ corresponds to the algorithm presented in~\cite{aloise2012improved}, where no aggregation of covering constraints is performed;
\item $\ell=n/2$ corresponds to the case where the initial aggregation $Q$ is obtained from a partition obtained by $k$-means using $\tilde{k} = n/2$;
\item $\ell=n/4$ corresponds to the case where the initial aggregation $Q$ is obtained from a partition obtained by $k$-means using $\tilde{k} = n/4$;
\item $\ell=k$ corresponds to the case where the initial aggregation $Q$ is made equal to the partition obtained by $k$-means for $k$ clusters, which is our maximum aggregation level.
\end{itemize}

We employ the same experimental setup for each CG+DCA execution. At each CG iteration, we add the minimum reduced cost compatible column, if available. Otherwise, we disaggregate the aggregating partition $Q$ using the incompatible column with the smallest reduced cost. We disaggregate the aggregated dual multipliers with the average strategy~\eqref{eq:disagg_avg}. Finally, a time limit of 7200 seconds was imposed for each CG+DCA execution.

Tables~1-4 of the Appendix C show that aggregating data points is always useful since the size of the aggregated RMP becomes smaller throughout the DCA process. Thus, the aggregation reduces considerably the proportion of the total time spent in the master problem resolution -- with a single exception between $\ell=n$ and $\ell=n/2$ for $k=2$ in data instance \texttt{gr202}.
Notably, the aggregation level $\ell=k$ impressively shows a time reduction in the range [55.58, 88.51]\% across all the values of $k$ for $\texttt{ch150}$, [55.30, 87.45]\% for $\texttt{gr202}$, [83.11, 91.73]\% for $\texttt{pr299}$, and [99.01, 99.79]\% for $\texttt{fl417}$. In particular, on $\texttt{pr299}$ with $k=2$ and $\texttt{fl417}$ with $k \in \{2, 3\}$ the baseline reaches the time limit of 2 hours, whereas our CG+DCA algorithm with different aggregation levels solves the problem very efficiently in a matter of seconds.  
Moreover, the number of iterations decreases considerably with respect to the baseline where DCA is not performed (aggregation level $\ell=n$), showing that the problem becomes less degenerate. 

We also observe that the final size of the aggregating partition $Q$ (column $m_{end}$) is often close to the original number of data points $n$. As a consequence, the column ``$Q$ updates'' shows more frequent disaggregation procedures when the initial number of components (column $m_{start}$) is smaller. However, it is noteworthy that this behavior of $m_{end}$ appears to vary depending on the data instance. For example, Table~4 shows significantly smaller values of $m_{end}$ as CG+DCA increases the aggregation level (column $\ell$) for instance \texttt{fl417}, which clearly results in the most substantial benefits of the DCA approach compared to the method previously proposed in~\cite{aloise2012improved}.

In view of the presented results, we can state that the column generation algorithm proposed in~\cite{aloise2012improved} benefits from the incorporation of the dynamic constraint aggregation technique of~\cite{elhallaoui2005dynamic, elhallaoui2010multi} at different aggregation levels. Therefore, we will use $\ell = k$ given its superiority performance observed in Tables~1-4 of Appendix C.


\subsection{DCA design choices}

In this section, we conduct ablation studies about possible design choices of the DCA methodology.
Ultimately, we are interested in knowing which of these choices leads to the best version of the CG+DCA algorithm for MSSC.
Our analysis is focused on evaluating  strategies for (i)
disaggregating the dual variables for the pricing problem at each iteration of the CG+DCA algorithm
(Section~\ref{sec:disaggregation}), and for (ii) 
selecting negative reduced cost incompatible columns for updating the aggregating partition $Q$ (Section~\ref{sec:q_update}).

\subsubsection{Dual variable disaggregation}

A procedure must be employed at each iteration of  CG+DCA  in order to disaggregate the dual solution obtained after solving the agRMP. Our first ablation study evaluates the use of three different strategies mentioned in Section~\ref{sec:disaggregation}, namely:
\begin{itemize}
\item \textit{sparse}: where a single dual variable associated with a component of $Q$ is randomly chosen to assume the entire value of the aggregated dual variable, while the other dual variables associated with that component are set to zero.

\item \textit{complementary}: that solves the LP complementary problem at each CG iteration for dual variable disaggregation (the reader is referred to~\cite{bouarab2017dynamic} for details).

\item \textit{average}: where the value of the aggregated dual variable is evenly distributed among the individual dual variables associated with that component.
\end{itemize}

 The bar charts presented in Online Appendix D show that the \textit{average} disaggregation strategy is better in terms of computing times since it does not require the solution of an additional linear program at every CG iteration, as opposed to the \textit{complementary} strategy.  Moreover, it appears to favor the generation of compatible columns, thus preventing unnecessary disaggregations of the aggregating partition $Q$ throughout the CG+DCA algorithm. Notably, the number of CG iterations and the number of constraints in the aggregated RMP are generally smaller, leading to substantial computational savings, with the only exception found when solving data instance \texttt{gr202} with $k=5$. In that case, $m_{avg}$  {(the average number of aggregated covering constraints per iteration)} and $m_{end}$  {(the final number of aggregated constraints in the final aggregated RMP)} are larger for the \textit{average} strategy than for the \textit{complementary} one. However, this did not lead to smaller computing times as the \textit{average} strategy yields less CG iterations besides dismissing the solution of additional LPs.

Regarding the \textit{sparse} strategy, the results reveal that this strategy, despite reducing the complexity of the pricing problem due to the sparsity of dual variables, does not yield satisfactory results. In particular, the pricing procedure fails to generate compatible columns, resulting in excessively rapid disaggregations of $Q$, consequently not leveraging DCA's advantages. The latter observation is confirmed by the statistics $m_{avg}$ and $m_{end}$, which become almost equal to the original size $n$. 

Our results indicate that the \textit{average} strategy is the most effective. Consequently, we will use it for the remainder of our computational experiments.
Finally, it is noteworthy that the value of $m_{end}$ in all our experiments never reached $n$ when employing the \textit{average} dual variable disaggregation strategy, thus revealing that several subsets of dual variables share the same value in the optimal dual solution. This finding certainly deserves future investigation.
According to our results, the ratio $m_{end}/n$ seems to depend on the characteristics of the data instance and the number of clusters sought.

\subsubsection{Aggregating partition disaggregation}
If no compatible columns with negative reduced cost are generated in the pricing procedure, then the partition $Q$ must be modified shall a negative reduce cost incompatible column exists.
To pivot an incompatible variable into the basis of
the aggregated problem, the aggregating partition must be modified by increasing its number of components and, thus, the number of constraints in the agRMP. 
Our second ablation study aims to evaluate two methodologies for selecting the candidate incompatible column for updating $Q$:
\begin{itemize}
    \item \textit{Min-RC}: adds to the RMP the minimum (negative) reduced cost incompatible column;
    \item \textit{Min-INC}: adds to the RMP the negative reduced cost incompatible column with the smallest number of incompatibilities (as described in Section \ref{sec:q_update}). Ties are broken by using the minimum reduced cost if necessary.  
\end{itemize}

The bar charts presented in Online Appendix E show that \textit{Min-INC} is more efficient than \textit{Min-RC}, primarily due to the smaller average number of constraints in the aggregated RMP (metric $m_{avg}$), leading to reduced computational times. This is due to the significantly smaller values of $u_{avg}$ {(average number of incompatibilities of the entering incompatible column)} when using \textit{Min-INC}, which leads to a smoother progression in the expansion of partition $Q$.
Conversely, using highly incompatible variables results in larger increases in the size of the aggregated RMP with each update, and consequently,  in a smaller number of updates (metric $Q$ updates) too.
This explains the particular exception when assessing \textit{Min-RC} and \textit{Min-INC} on  dataset \texttt{gr202} with  
$k=5$ clusters. In that case, despite $m_{avg}$ being larger with the \textit{Min-RC} strategy compared to the \textit{Min-INC} strategy, the CPU times are slightly larger for the latter. This increase is attributed to the time spent on $Q$ updates (and on updating the corresponding agRMP), which are not offset in this particular instance by smaller agRMPs on average (118.15 vs. 123.43) and slightly lower number of CG iterations (2813 vs. 2836).
Given the superiority of \textit{Min-INC} over \textit{Min-RC} according to our experiments, we will use \textit{Min-INC} for the final benchmarked version of CG+DCA described in the next section.

\subsection{Comparison with state-of-the-art solver}

This section aims to compare our CG+DCA method for MSSC against the current state-of-the-art solver SOS-SDP \citep{piccialli2022sos} found in the literature. SOS-SDP is an exact algorithm based on the branch-and-cut technique that features semidefinite programming bounds. For SOS-SDP, we consider all the parameter settings recommended
in \cite{piccialli2022sos} and the code, written in C++, is available online\footnote{\url{https://github.com/antoniosudoso/sos-sdp}}.
Regarding our CG+DCA algorithm, at each CG iteration, at most 10 compatible columns with negative reduced costs 
are added to the RMP. Following our ablation studies,  the dual variables values are disaggregated by using the \textit{average} strategy, the initial aggregation level of partition $Q$ is set to $\ell = k$, and the partition $Q$ is updated using the \textit{Min-INC} strategy. When branching is necessary, the enumeration tree is visited with the best-first search strategy. Each child node inherits from its parent the aggregating partition $Q$ and the set of columns in the optimal basis
that do not violate the branching constraint. Additionally, we augment this initial set of columns with positive reduced cost columns whose values fall within the integrality gap (see, e.g.~\cite{baldacci2008exact,ronnberg2019integer}) up to a limit of $2m $ columns.

The dataset instances used in  tests of this section are divided into two groups: 
\begin{itemize}
\item \textit{small to medium-size instances}: \texttt{pr299},
\texttt{fl417},
\texttt{ali535}, and
\texttt{gr666};

\item \textit{large-scale instances}: 
\texttt{u1060},
\texttt{u2152},
\texttt{fl3795}, and
\texttt{rl5934}.
\end{itemize}

For the small and medium-size instances we tested both algorithms CG+DCA and SOS-SDP for $k \in \{2, 4, 6, 8, 10\}$. The large-scale instances are solved for $k \in \{r, r+10, r+20, r+30, r+40\}$, where $r$ is the multiple of 10 closest to $\lfloor n/100 \rfloor$. Moreover, the exact solution of the small and medium-size instances considers an
optimality tolerance on the percentage gap of $\varepsilon = 0.01\%$. That is, we terminate the branch-and-bound method when $(\textrm{UB} - \textrm{LB}) / \textrm{UB} \leq \varepsilon$, where UB and LB denote the best upper and lower bounds, respectively.
For the large-scale instances, the optimality tolerance 
$\varepsilon$ is set to $0.1\%$ in accordance with the value used by the SOS-SDP solver in~\cite{piccialli2022sos}. Moreover, we set a time limit of 86,400 seconds (24 hours) of CPU time.

\begin{table}[!ht]
    \centering
    \footnotesize
     \caption{Comparison between SOS-SDP and CG+DCA on small and medium-size instances.}
\begin{tabular*}{\textwidth}{@{\extracolsep{\fill}}lccllllll@{\extracolsep{\fill}}}
\toprule
& & & \multicolumn{3}{c}{SOS-SDP} & \multicolumn{3}{c}{CG+DCA}\\
\cmidrule{4-6}\cmidrule{7-9}%
Data & $k$	&	$f_{\textrm{opt}}$	&	$\textrm{Gap}_0 (\%)$	&	Nodes	&	Time (s) &	$\textrm{Gap}_0 (\%)$	&	Nodes	&	Time (s) \\
\midrule
\multirow{6}{*}{$\texttt{pr299}$}	&	2	&	4.00724E+08	&	0.15	&	5	&	68.69	&	0.01	&	1	&	60.74	\\
	&	4	&	2.17262E+08	&	1.19	&	73	&	906.77	&	0.00	&	1	&	88.73	\\
	&	6	&	1.35426E+08	&	1.31	&	97	&	1255.27	&	0.01	&	1	&	29.42	\\
	&	8	&	9.93752E+07	&	0.55	&	35	&	429.58	&	0.73	&	5	&	127.75	\\
	&	10	&	7.33670E+07	&	0.84	&	21	&	338.12	&	0.00	&	1	&	21.43	\\ 
 \midrule	
\multirow{6}{*}{$\texttt{fl417}$}	&	2	&	1.07735E+08	&	0.00	&	1	&	22.98	&	0.00	&	1	&	49.34	\\
	&	4	&	3.66438E+07	&	6.04	&	57	&	1432.39	&	0.00	&	1	&	8.36	\\
	&	6	&	1.29071E+07	&	0.97	&	13	&	243.12	&	0.00	&	1	&	8.13	\\
	&	8	&	7.62489E+06	&	1.08	&	17	&	334.84	&	0.00	&	1	&	11.39	\\
	&	10	&	5.53184E+06	&	0.00	&	1	&	216.93	&	0.00	&	1	&	10.31	\\ 
  \midrule
\multirow{6}{*}{$\texttt{ali535}$}	&	2	&	9.90552E+05	&	0.39	&	15	&	977.48	&	0.00	&	1	&	849.72	\\
	&	4	&	3.20279E+05	&	1.01	&	43	&	1243.10	&	0.00	&	1	&	192.75	\\
	&	6	&	2.00494E+05	&	0.80	&	39	&	1060.07	&	0.00	&	1	&	744.23	\\
	&	8	&	1.47785E+05	&	5.13	&	179	&	4896.03	&	0.10	&	3	&	739.30	\\
	&	10	&	1.09825E+05	&	4.65	&	163	&	4558.82	&	0.00	&	1	&	573.26	\\
 \midrule
\multirow{6}{*}{$\texttt{gr666}$}	&	2	&	1.75401E+06	&	5.44	&	197	&	10423.47	&	0.00	&	1	&	2904.74	\\
	&	4	&	6.13995E+05	&	3.09	&	111	&	7356.33	&	0.01	&	1	&	3516.18	\\
	&	6	&	3.82677E+05	&	2.00	&	81	&	6218.04	&	0.01	&	1	&	1893.98	\\
	&	8	&	2.85925E+05	&	4.47	&	165	&	9385.76	&	0.01	&	1	&	1006.45	\\
	&	10	&	2.24184E+05	&	7.11	&	257	&	14553.21	&	0.01	&	1	&	963.42	\\
\bottomrule
    \end{tabular*}
    \label{tab:result_cg_vs_sdp_small}
\end{table}

In Tables \ref{tab:result_cg_vs_sdp_small} and \ref{tab:result_cg_vs_sdp_large} we present computational results of the two methods, SOS-SDP and CG+DCA, across the datasets with varying numbers of clusters. For each dataset and cluster number ($k$), we report the optimal MSSC solution value ($f_{\textrm{opt}}$), the percentage gap at the root node ($\textrm{Gap}_0 (\%)$), the number of explored nodes (Nodes), and the computational times (Time).  

The results in Table \ref{tab:result_cg_vs_sdp_small}  demonstrate that our CG+DCA algorithm consistently outperforms SOS-SDP, particularly in larger datasets where the computational time of SOS-SDP increases significantly. For example, in the $\texttt{gr666}$ dataset with $k=2$, CG+DCA closes the gap in 2904.74 seconds, while SOS-SDP requires 10423.47 seconds, indicating a speed advantage of 3.6 times. For $k=4$ clusters, CG+DCA finishes in 3516.18 seconds, compared to 7356.33 seconds for SOS-SDP, making CG+DCA approximately twice as fast. The $\texttt{ali535}$ dataset further highlights CG+DCA's efficiency; for $k=8$, CG+DCA closes the gap in 739.30 seconds, whereas SOS-SDP takes 4896.03 seconds, making CG+DCA about seven times faster. With $k=10$, CG+DCA solves the problem in 573.26 seconds, significantly less than the 4558.82 seconds required by SOS-SDP, showcasing an eight-fold improvement in speed. The only exception occurs in the $\texttt{fl417}$ dataset for $k=2$, where both algorithms solve the problem at the root node, but SOS-SDP performs slightly better, closing the gap in 22.98 seconds compared to 49.34 seconds for CG+DCA. Although the number of nodes explored by the two methods is not directly comparable due to their distinct bounding nature, SOS-SDP tends to show higher gaps at the root node and thus explores more nodes. Notably, for CG+DCA, the root node gap is smaller than 0.01\% in most instances, indicating that the algorithm often terminates directly at the root without the need for branching. This confirms that when the CG+DCA algorithm converges, the solution of the LP relaxation is often integer and thus optimal. In only two instances, namely
$\texttt{pr299}$ for $k=8$, and $\texttt{ali535}$ for $k=8$, a small gap exists at the root node, and thus branching becomes necessary. Interestingly, only a few nodes are explored in these cases, highlighting the effectiveness of the branching strategy.

Finally, we remark that due to the use of DCA, the aggregated RMP is typically less degenerate and has a significantly smaller size compared to the non-aggregated RMP. This allows the addition of multiple columns without compromising the efficiency of the solution process. In fact, adding multiple columns with negative reduced costs at each CG iteration reduces the number of iterations, thereby decreasing the CPU time required to obtain an optimal solution. Comparing the computational results in Tables 3-4 in Online Appendix C for the  $\texttt{pr299}$ and $\texttt{fl417}$ datasets, respectively, with aggregation level $\ell = k$, to those in Table \ref{tab:result_cg_vs_sdp_small}, it is evident that the computational time decreases substantially. Similar observations apply to the $\texttt{ali535}$ and $\texttt{gr666}$ datasets.

In Table \ref{tab:result_cg_vs_sdp_large} we present a comparison between SOS-SDP and CG+DCA on large-scale instances. 
When the 24-hour time limit is reached, the ``Nodes'' column shows the achieved percentage gap, indicated in parenthesis.
Computational results indicate that the CG+DCA method effectively solves all instances within the desired optimality gap. Similarly to the results reported in Table \ref{tab:result_cg_vs_sdp_small}, the majority of instances are solved directly at the root node. Only in two cases, for $\texttt{u2172}$ with $k=50$ and $\texttt{fl3795}$ with $k=40$, a gap exists at the root node, necessitating branching. However, even in these cases, only three nodes are explored before the gap is closed. In contrast, the SOS-SDP method demonstrates mixed results. While it successfully closes the gap for all instances in the $\texttt{u1060}$ dataset and for $k \in \{20, 30, 40\}$ in the $\texttt{u2172}$ dataset, it requires significantly more CPU time compared to CG+DCA. For the remaining instances in the $\texttt{u2172}$ dataset and for all $k$ values in the $\texttt{rl5934}$ dataset, SOS-SDP reaches the time limit with substantial gaps remaining. Moreover, for the $\texttt{fl3795}$ dataset SOS-SDP explores only a minimal number of nodes - and only~1 for the $\texttt{rl5934}$ dataset - before the time limit is reached. We highlight that the poor performance of SOS-SDP on large-scale instances is not unexpected. Although this solver has successfully handled instances with a few thousand data points, its bottleneck lies in the resolution of SDP relaxations, which are further complicated by the addition of valid inequalities \citep{piccialli2022sos}. 

Overall, the results reported in this section highlight the efficiency and robustness of the CG+DCA approach for large-scale instances, showing improvements in computational time ranging from approximately 3 times faster (e.g., $\texttt{u1060}$, $k=10$) to more than 10 times faster (e.g., $\texttt{fl3795}$, $k=80$).

\begin{table}[!ht]
    \centering
    \footnotesize
    \caption{Comparison between SOS-SDP and CG+DCA on large-scale instances.}
\begin{tabular*}{\textwidth}{@{\extracolsep{\fill}}lcclllccl@{\extracolsep{\fill}}}
\toprule
& & & \multicolumn{3}{c}{SOS-SDP} & \multicolumn{3}{c}{CG+DCA}\\
\cmidrule{4-6}\cmidrule{7-9}%
data & $k$	&	$f_{\textrm{opt}}$	&	$\textrm{Gap}_0 (\%)$	&	Nodes	&	Time (h) &	$\textrm{Gap}_0 (\%)$	&	Nodes	&	Time (h) \\
\midrule
\multirow{5}{*}{$\texttt{u1060}$}	&	10	&	1.75484E+09	&	6.08	&	153	&	6.53	&	0.01	&	1	&	1.87	\\
	&	20	&	7.92038E+08	&	4.10	&	97	&	3.77	&	0.06	&	1	&	0.23	\\
	&	30	&	4.81252E+08	&	5.09	&	129	&	4.44	&	0.00	&	1	&	0.09	\\
	&	40	&	3.41343E+08	&	5.48	&	141	&	5.96	&	0.04	&	1	&	0.06	\\
	&	50	&	2.55510E+08	&	5.32	&	155	&	6.15	&	0.03	&	1	&	0.04	\\ \midrule
\multirow{5}{*}{$\texttt{u2172}$}	&	20	&	7.45883E+07	&	4.61	&	55	&	18.15	&	0.10	&	1	&	10.29	\\
	&	30	&	4.64093E+07	&	4.76	&	63	&	17.64	&	0.01	&	1	&	6.51	\\
	&	40	&	3.36536E+07	&	6.35	&	91	&	22.03	&	0.09	&	1	&	4.3	\\
	&	50	&	2.56488E+07	&	8.11	&	75 (4.03 \%)	&	24.00	&	0.14	&	3	&	2.33	\\
	&	60	&	2.08484E+07	&	8.74	&	67 (6.77 \%)	&	24.00	&	0.07	&	1	&	0.44	\\ \midrule
\multirow{5}{*}{$\texttt{fl3795}$}	&	40	&	1.27082E+07	&	7.68	&	4 (5.99 \%)	&	24.00	&	0.33	&	3	&	16.79	\\
	&	50	&	8.84954E+06	&	14.74	&	7 (12.43 \%)	&	24.00	&	0.09	&	1	&	5.63	\\
	&	60	&	6.45195E+06	&	13.62	&	5 (10.38 \%)	&	24.00	&	0.08	&	1	&	3.77	\\
	&	70	&	4.98172E+06	&	17.25	&	3 (15.91 \%)	&	24.00	&	0.01	&	1	&	2.68	\\
	&	80	&	4.09405E+06	&	17.98	&	9 (16.44 \%)	&	24.00	&	0.01	&	1	&	2.24	\\ \midrule
\multirow{5}{*}{$\texttt{rl5934}$}	&	60	&	2.62478E+09	&	9.36	&	1 (9.36 \%)	&	24.00	&	0.07	&	1	&	23.36	\\
	&	70	&	2.20020E+09	&	12.49	&	1 (12.49 \%)	&	24.00	&	0.10	&	1	&	22.12	\\
	&	80	&	1.90220E+09	&	11.22	&	1 (11.22 \%)	&	24.00	&	0.05	&	1	&	20.09	\\
	&	90	&	1.66614E+09	&	18.36	&	1 (18.36 \%)	&	24.00	&	0.10	&	1	&	18.97	\\
	&	100	&	1.47966E+09	&	16.78	&	1 (16.78 \%)	&	24.00	&	0.09	&	1	&	17.98	\\
\bottomrule
    \end{tabular*}
    
    \label{tab:result_cg_vs_sdp_large}
\end{table}

\section{Conclusions}\label{sec:final}

In this paper, we have proposed a CG algorithm enhanced by DCA for solving the MSSC problem. 
In particular, our focus was on accelerating the resolution of the RMP introduced in~\cite{aloise2012improved}, which is known to suffer from high degeneracy. 
Our computational experiments highlight that incorporating DCA effectively mitigates the impact of degeneracy when solving MSSC instances using branch-and-price, enabling the resolution of data instances in the plane with up to 5934 data points -- a new record for MSSC global optimization. Remarkably, our approach was found to be 3 to 10 times faster than the current state-of-the-art exact method on the tested instances.

Our algorithm for MSSC was fine-tuned through a series of ablation studies, which demonstrated the superiority of (i) disaggregating dual variables by equally distributing their values among the corresponding aggregated dual variables, and (ii) by adding negative reduced cost columns with minimum number of incompatibilities with respect to the current constraint aggregation partition. These studies are important to consolidate the use of DCA in broader mathematical optimization contexts, allowing the identification of its strengths and possible improvements. Indeed, other aspects of DCA are still to be investigated such as eventually downsizing the partition of aggregated constraints during the execution of the CG algorithm.

Finally, as future work, we aim to extend the proposed CG+DCA methodology to solve general multidimensional MSSC instances arising in data mining and unsupervised machine learning applications.

\section*{Acknowledgments}
First and foremost, we express our deepest gratitude to Professor François Soumis (1946–2025) for his invaluable insights throughout the development of this research. The authors also gratefully acknowledge the financial support of the Natural Sciences and Engineering Research Council of Canada (NSERC).

\bibliography{sn-bibliography}

\end{document}